\documentclass{amsart}
\usepackage{amsmath,amscd,amssymb,amsthm,amsbsy,amscd,stmaryrd}
\usepackage[dvips]{graphicx,color}
\usepackage{epsfig}
\usepackage{mathrsfs}
\usepackage{mathpazo}
\input{xy}

\newtheorem{theorem}{Theorem}[section]
\newtheorem{lemma}[theorem]{Lemma}
\newtheorem{proposition}[theorem]{Proposition}

\theoremstyle{definition}

\newtheorem{remark}[theorem]{Remark}

\numberwithin{equation}{section}

\begin{document}
\title[Local curvature estimate]{A local curvature estimate for the Ricci-harmonic flow on complete Riemannian manifolds}
\author{Yi Li}
\address{School of Mathematics and Shing-Tung Yau Center, Southeast University, Nanjing 211189, China}
\email{yilicms@gmail.com; yilicms@seu.edu.cn}

\author{Miaosen Zhang}
\address{Chien-Shiung Wu College of SEU, Southeast University, Nanjing 211189, China}
\email{zhmise@163.com}



\maketitle


\begin{abstract}
In this paper we consider the local $L^p$ estimate of Riemannian curvature for the Ricci-harmonic flow or List's flow introduced by List \cite{List2005} on complete noncompact manifolds. As an application, under the assumption that the flow exists on a finite time interval $[0,T)$ and the Ricci curvature is uniformly bounded, we prove that the $L^p$ norm of Riemannian curvature is bounded, and then, applying the De Giorgi-Nash-Moser iteration method, obtain the local boundedness of Riemannian curvature and consequently the flow can be continuously extended past $T$.
\end{abstract}

\section{Introduction}\label{section1}

The Ricci-harmonic flow is defined to be the following system:
\begin{align}
    \begin{cases}
    \displaystyle{\frac{\partial }{\partial t}g(t) = -2\!\ {\rm Ric}(g(t)) +4\!\ du(t)\otimes du(t)},\\
    \\
     \displaystyle{\frac{\partial }{\partial t}u(t) = \Delta_{g(t)} u(t)},\\
     \\
    g(0) = g_0,\quad u(0) = u_0,
    \end{cases},
\end{align}
where $g_{0}$ is a fixed Riemannian metric, $u_{0}$ is a fixed smooth function, $t\in [0,T)$, $g(t)$ is a family of metrics, $u=u(t)$ is a family of smooth functions on an $n$-dimensional manifold $M$. It was first introduced in \cite{List2005} and also called extended Ricci flow in \cite{A2020, HL2018, List2005, WZ2019}. The flow equations, as the motivation for studying it, were proved to characterize the static Einstein vacuum metrics \cite{EK1962, List2005}. Under the assumption that $M$ is compact, List \cite{List2005, List2008} prove the short time existence, and also proved that if the Riemann curvature is uniformly bounded for all $t\in [0,T)$, then the solution
can be extended beyond $T$.  For a more general setting, see \cite{Muller, Muller2012}. In the complete noncompact case, the long time existence of manifolds with bounded scalar curvature was given by the first author \cite{Li2018}.

\par Over the last decade, there are lots of works on both compact and noncompact manifolds about eigenvalues, entropies, functionals, and solitons, see, for example, \cite{AAB2019, AKO2019, A2020, CGT2015, FZ2016, GPA2013, GPA2014a, GPA2014b, HL2018, LY2014, LW2017, WZ2019, YS2012}. In this paper, we mainly focus on the estimate of curvature. List \cite{List2008} proved that, $M$ being compact, the Ricci-harmonic flow can be extended if the Riemannian curvature is bounded, as an application to see the importance of curvature estimate. Unfortunately, counterexamples show that the Riemannian curvature (see \cite{List2005}) and Ricci curvature (see \cite{CZ2013}) could not be bounded without any restrictions. On the other hand, those curvatures are $L^2$ bounded in certain cases (e.g. $n=4$) if scalar curvature is bounded (see \cite{LY2018}). Furthermore, the pseudo-locality theorem corresponding to the Ricci-harmonic flow was be given in \cite{GHP2018}. However, as in the Ricci flow case, whether the scalar curvature is bounded for the Ricci-harmonic flow remains an open problem (see in \cite{LY2016}).

Instead of giving a point-wise estimate of ${\rm Rm}$, the $L^p$ norm
\begin{align*}
    ||{\rm Rm}||_{p,M\times [0,T)} = \left(\int_0^T\int_M|{\rm Rm}(g(t))|^p_{g(t)}dV_{g(t)}dt\right)^{\frac{1}{p}}
\end{align*}
for ${\rm Rm}$ was recently established in \cite{WZ2020} on compact manifolds. The main result of this paper is to give a local $L^p$ and point-wise estimate for complete manifolds, strengthening the propositions in \cite{Li2018}.

${}$

\textbf{Notations:} In the following, we often omit $t$ variable, for example, $g=g(t)$, $u=u(t)$, $\Delta=\Delta_{g(t)}$, etc. The operator $\square := \partial_t - \Delta$ will be frequently used later. $C$ represents positive finite constants that we don't care about their value.

${}$

The first result of this paper is

\begin{theorem}\label{t1.1} {\rm (also see Theorem \ref{t2.7})} Let $(g(t),u(t))_{t\in[0,T]}$ be a solution to the Ricci-harmonic flow on $M\times[0,T]$, where $M$ is a complete $n$-dimensional manifold and $T\in(0,+\infty)$. Suppose there exist constants $\rho,K,L>0$ and a point $x_0\in M$ such that the geodesic ball $B_{g(0)}(x_0,\rho/\sqrt{K})$ is compactly contained on $M$ and
\begin{align}
|{\rm Ric}(g(t))|_{g(t)} \leq K ,\qquad |\nabla_{g(t)}u(t)|_{g(t)}\leq L.\label{1.2}
\end{align}
 For any $p\geq 3$, there exist constants $\Gamma_1,\Gamma_2$ depending only on $n,p,\rho,K$, $L$ and $T$, such that
\begin{align*}
\int_{B_{g(0)}(x_0,\rho/2\sqrt{K})}|{\rm Rm}(g(t))|_{g(t)}^p&dV_{g(t)}\leq \Gamma_1\int_{B_{g(0)}(x_0,\rho/\sqrt{K})}|{\rm Rm}(g(0))|^p_{g(0)}dV_{g(0)}\\&+\Gamma_2{\rm Vol}_{g(0)}\left(B_{g(0)}\left(x_0,\frac{\rho}{\sqrt{K}}\right)\right).
\end{align*}
\end{theorem}

Actually the explicit expressions for $\Gamma_{1}$ and $\Gamma_{2}$ can be found in the proof of Theorem \ref{t2.7}.

${}$

Under the additional condition that $|\nabla^2_{g(t)}u|_{g(t)}$ is bounded, Theorem 1.1 was proved in \cite{Li2018}. Theorem \ref{t1.1} shows that this additional condition can be removed. According to the following remark, the boundedness of $|\nabla_{g(t)} u(t)|_{g(t)}$ can also be removed. We include the condition $|\nabla_{g(t)}u(t)|_{g(t)}\leq L$ in
Theorem \ref{t1.1} is in order to see how $K$ and $L$ involve in the $L^p$ estimate of ${\rm Rm}$.

\begin{remark}\label{r1.2} {\rm (see Theorem B.2 in \cite{Li2018})} Suppose that $(g(t),u(t))_{t\in[0,T]}$ is a solution to (1) on $M\times [0,T]$, where $M$ is a complete $n$-dimensional manifold. If the estimate
\begin{align*}
    \underset{M\times [0,T)}{\sup} |{\rm Ric}(g(t))|_{g(t)} \leq K
\end{align*}
holds for some positive constant $K$, then we have
\begin{align*}
    \underset{M\times [0,T)}{\sup} |\nabla_{g(t)} u(t)|^2_{g(t)} \leq 2KC(n),
\end{align*}
where $C(n)$ is a positive number depends only on $n$.
\end{remark}

Theorem \ref{t1.1} and Remark \ref{r1.2} imply

${}$

\begin{theorem}\label{t1.3} Let $(g(t),u(t))_{t\in[0,T]}$ be a solution to the Ricci-harmonic flow on $M\times[0,T]$, where $M$ is a complete $n$-dimensional manifold and $T\in(0,+\infty)$. Suppose there exist constants $\rho,K$ and a point $x_0\in M$ such that the geodesic ball $B_{g(0)}(x_0,\rho/\sqrt{K})$ is compactly contained on $M$ and
\begin{align}
|{\rm Ric}(g(t))|_{g(t)} \leq K.\label{1.3}
\end{align}
For any $p\geq 3$, there exist constants $\Gamma_1,\Gamma_2$ depending only on $n,p,\rho, K$, and $T$ such that
\begin{eqnarray*}
\int_{B_{g(0)}(x_0,\rho/2\sqrt{K})}|{\rm Rm}(g(t))|_{g(t)}^p\!\ dV_{g(t)}&\leq& \Gamma_1\int_{B_{g(0)}(x_0,\rho/\sqrt{K})}|{\rm Rm}(g(0))|^p_{g(0)}dV_{g(0)}\\
&&+ \ \Gamma_2{\rm Vol}_{g(0)}\left(B_{g(0)}\left(x_0,\frac{\rho}{\sqrt{K}}\right)\right).
\end{eqnarray*}
\end{theorem}

Finally we state our main theorem.

\begin{theorem}\label{t1.4} {\rm (also see Theorem 3.2)} Let $(g(t),u(t))_{t\in[0,T)}$ be a smooth solution to the Ricci-harmonic flow on $M\times[0,T)$ with $T\in(0,+\infty)$, where $M$ is a complete $n$-dimensional manifold. If $(M,g(0))$ is complete and:
    \begin{align*}
        \underset{\quad M\quad }{\sup}|{\rm Rm}(g(0))|_{g(0)}<\infty, \ \ \  \sup_{M\times[0,T)}|{\rm Ric}(g(t))|_{g(t)}<\infty
    \end{align*}
    then the flow can be extended over $T$.
\end{theorem}

This paper is organized as follow: In Sect. 2.1, we state our main idea and prove Theorem \ref{t1.1}, i.e., the $L^p$ norm estimate of Riemannian curvature. We supply the details of the proof in Sect. 2.2. In Sect. 3, We discuss the extension of (1.1) and prove Theorem \ref{t1.4}.

\section{$L^p$ estimate of Riemannian curvature}\label{section2}

We start with the proof of Theorem \ref{t1.1}. As in \cite{KMW2016, Li2018}, we let $\phi$ be a (time independent) Lipschitz function with compact support in a domain $\Omega\subset M$. Throughout this section, we always assume the condition (\ref{1.2}) holds.

\subsection{Main idea}\label{subsection2.1}

Given a real number $p\geq1$ that is determined later. We introduce the following integrals:
\[B_1 := \frac{1}{K}\int_M|\nabla {\rm Ric}|^2|{\rm Rm}|^{p-1}\phi^{2p}dV_t,\quad B_2 := \int_M|\nabla {\rm Rm}|^2|{\rm Rm}|^{p-3}\phi^{2p}dV_t,\]
and also
\[A_1 := \int_M|{\rm Rm}|^{p}\phi^{2p}dV_t,\quad A_2 := \int_M|{\rm Rm}|^{p-1}\phi^{2p}dV_t,\]
\[A_3 := \int_M|{\rm Rm}|^{p-1}|\nabla\phi|^2\phi^{2p-1}dV_t, \quad A_4:= \int_M|{\rm Rm}|^{p-1}|\nabla\phi|^2\phi^{2p-2}dV_t.\]
In order to control the second derivative of $u$, we
need another type of integrals
\[T_k := \int_M|{\rm Rm}|^{k-1}|\nabla^2u|^2\phi^{2p}dV_t, \qquad k=1,2, \cdots, p.\]

Then we have following inequalities, proved in Sect. \ref{subsection2.2}.

\begin{proposition}\label{p2.1} We have
\[\frac{d}{dt}A_{1}\leq B_1+CKB_2+CKA_4+C(K+L^2)A_1+CT_p\]
\end{proposition}

\begin{proposition}\label{p2.2}
\begin{eqnarray*}
    B_1&\leq& CKB_2+C(K+L^2)A_1+CKL^2A_2+CKA_4\\
    &&+ \ CT_p-\frac{1}{2K}\frac{d}{dt}\left(\int_M|{\rm Ric}|^2|{\rm Rm}|^{p-1}\phi^{2p}dV_t\right).
\end{eqnarray*}
\end{proposition}

We observe that all $T_{k}$ can be controlled by $T_{p}$ and $T_{1}$.

${}$

\begin{lemma}\label{l2.3} For any positive constant C and any $k = 1,2,...,p$,
\[T_k\leq \frac{1}{C^{p-k}}T_p+(p-k)C^{k-1}T_1\]
\end{lemma}

\begin{proof} We can easily find that, for any positive constant $C$, the following inequality
\[(|{\rm Rm}|-C)(|{\rm Rm}|^{k-1}-C^{k-1})\geq 0,\]
holds, which implies
\[|{\rm Rm}|^{k}-C|{\rm Rm}|^{k-1}+C^k\geq C^{k-1}|{\rm Rm}|\geq 0.\]
Integrating on both sides yields
\[T_k\leq \int_M\left(\frac{1}{C}|{\rm Rm}|^{k}+C^{k-1}\right)|\nabla^2u|^2\phi^{2p}dV_t=\frac{1}{C}T_{k+1}+C^{k-1}T_1.\]
We now use the induction method to prove this lemma. For $k = p$, $T_p\leq T_p$ satisfied. If the lemma is satisfied for some $k\leq p$, then
\begin{align*}
    T_{k-1}&\leq \frac{1}{C}T_k+C^{k-2}T_1\\
    &\leq\frac{1}{C}\left(\frac{1}{C^{p-k}}T_p+(p-k)C^{k-1}T_1\right)+C^{k-2}T_1\\
    &=\frac{1}{C^{p-(k-1)}}T_p+[p-(k-1)]C^{k-2}T_1.
\end{align*}
Therefore the above mentioned estimate hold.
\end{proof}

According to Lemma \ref{l2.3}, we can estimate all $T_k$'s in terms of $T_p$
and $T_1$. However, from the definition, we see that $T_{p}$ and $T_{1}$ contain the second derivative of $u$ so that we can not use the condition (\ref{1.2}) to bound them. More precisely,
$$
T_{p}=\int_{M}|{\rm Rm}|^{p-1}|\nabla^{2}u|^{2}\phi^{2p}dV_{t}, \ \ \ T_{1}=\int_{M}|\nabla^{2}u|^{2}\phi^{2p}dV_{t}.
$$
Motivated by these two integrals, by replacing the second derivative of $u$ by its first derivative, we set
\begin{align*}
    S := \int_M|{\rm Rm}|^{p-1}|\nabla u|^2\phi^{2p}dV_t,\qquad \widetilde{S}=\int_M|\nabla u|^2\phi^{2p}dV_t,
\end{align*}
It is clear from the condition (\ref{1.2}) that $S\leq L^2A_2$.

${}$

\begin{proposition}\label{p2.4} We have
\[B_2\leq -\frac{1}{p-1}\frac{d}{dt}A_{2}+CA_1+CA_4+CL^2A_2+CT_{p-1}\]
\end{proposition}

\begin{proposition}\label{p2.5} For each $p\geq2$, $T_p$ satisfies the following estimate
\begin{align*}
    T_p\leq-\frac{d}{dt}S-&\frac{C}{p-1}\frac{d}{dt}A_2-\frac{C}{K}\frac{d}{dt}\left(\int_M|{\rm Ric}|^2|{\rm Rm}|^{p-1}\phi^{2p}dV_t\right)\\&+C(K+L^2)A_1+CKL^2A_2+C(K+L^2)A_4+C^{p-1}T_1
\end{align*}
\end{proposition}

\begin{proposition}\label{p2.6} $T_1$ satisfies the following estimate
\begin{align*}
    T_1\leq -\frac{d}{dt}\widetilde{S}+CL^2{\rm Vol}_{g(t)}(\Omega).
\end{align*}
\end{proposition}

We will give proofs for Proposition \ref{p2.4} -- Proposition \ref{p2.6} in Sect. \ref{subsection2.2}. Now we can prove Theorem \ref{t1.1}.

\begin{theorem}\label{t2.7} Let $(g(t),u(t))_{t\in[0,T]}$ be a solution to the Ricci-harmonic flow on $M\times[0,T]$, where $M$ is a complete $n$-dimensional manifold with $T\in(0,+\infty)$. Suppose that there exist constants $\rho,K,L>0$ and a point $x_0\in M$ such that the geodesic ball $B_{g(0)}(x_0,\rho/\sqrt{K})$ is compactly contained on $M$ and $({\rm Ric}(g(t)), \nabla u(t))$ satisfies (\ref{1.2}). For any $p\geq 3$, there exist constants $\Gamma_1,\Gamma_2$ depending only on $n,p,\rho,K$, $L$ and $T$, such that
\begin{align*}
\int_{B_{g(0)}(x_0,\rho/2\sqrt{K})}|{\rm Rm}(g(t))|_{g(t)}^p&dV_{g(t)}\leq \Gamma_1\int_{B_{g(0)}(x_0,\rho/\sqrt{K})}|{\rm Rm}(g(0))|^p_{g(0)}dV_{g(0)}\\&+\Gamma_2{\rm Vol}_{g(0)}\left(B_{g(0)}\left(x_0,\frac{\rho}{\sqrt{K}}\right)\right).
\end{align*}
\end{theorem}

Actually the explicit expressions for $\Gamma_{1}$ and $\Gamma_{2}$ can be found in the proof.

\begin{proof} Applying Lemma \ref{l2.3} with $C=1$ and $k=p-1$ to Proposition \ref{p2.4} yields
\begin{align*}
    B_2\leq -\frac{1}{p-1}\frac{d}{dt}A_{2}+CA_1+CA_4+CL^2A_2+CT_{p}+CT_1
\end{align*}
Plugging Proposition \ref{p2.2}, the above inequality into Proposition \ref{p2.1} successively to replace $B_1$ and $B_2$:
\begin{align*}
    \frac{d}{dt}&\left[A_{1}+\frac{CK}{p-1}A_{2}+\frac{1}{2K}\int_M|{\rm Ric}|^2|{\rm Rm}|^{p-1}\phi^{2p}dV_t\right]\\ &\qquad\qquad \leq C(K+L^2)A_1+CKL^2A_2+CKA_4+CKT_{p}+CKT_1
\end{align*}
Then apply proposition \ref{p2.5} and Proposition \ref{p2.6} to replace $T_{p}$ and $T_{1}$, we obtain
\begin{align*}
\frac{d}{dt}&\left[A_1+\frac{CK}{p-1}A_2+KC^{p}\widetilde{S}+CKS+C\int_M|{\rm Ric}|^2|{\rm Rm}|^{p-1}\phi^{2p}dV_t\right]\\& \qquad\qquad\leq CK(K+L^2)A_1+CK^2L^2A_2+CK(K+L^2)A_4+CKC^{p}L^2{\rm Vol}_{g(t)}(\Omega).
\end{align*}
Choose $\Omega:=B_{g(0)}\left(x_0,\rho/\sqrt{K}\right)$ and
\[\phi:=\left(\frac{\rho/\sqrt{K}-d_{g(0)}(x_0,\cdot)}{\rho/\sqrt{K}}\right)_+.\]
Define
\begin{align*}
U := \int_M|{\rm Rm}|^{p}\phi^{2p}dV_t&+\frac{CK}{p-1}\int_M|{\rm Rm}|^{p-1}\phi^{2p}dV_t +C\int_M| {\rm Ric}|^2|{\rm Rm}|^{p-1}\phi^{2p}dV_t \\&+CK\int_M|{\rm Rm}|^{p-1}|\nabla u|^2\phi^{2p}dV_t +KC^p\int_M|\nabla u|^2\phi^{2p}dV_t
\end{align*}
then $U$ satisfies the following estimate
\[U'\leq \left[CK^2+CKL^2+C(p-1)KL^2\right]U+CK(K+L^2)A_4+CKC^{p}L^2{\rm Vol}_{g(t)}(\Omega).\]
using
$$
e^{-2Kt}g(0)\leq g(t)\leq e^{2Kt}g(0)
$$
and
$$
|\nabla_{g(t)}\phi|_{g(t)}\leq e^{KT}|\nabla_{g(0)}\phi|_{g(0)}\leq \sqrt{K}e^{KT}/\rho.
$$
we can estimate $A_{4}$ as follows:
\begin{align*}
    A_4 = \int_M|{\rm Rm}|^{p-1}&|\nabla\phi|^2\phi^{2p-2}dV_t\leq\int_{B_{g(0)}(x_0,\rho/\sqrt{K})}|{\rm Rm}|^{p-1}\phi^{2p-2}K\rho^{-2}e^{2KT}dV_t\\&\leq\int_{B_{g(0)}(x_0,\rho/\sqrt{K})}\left[\frac{(|{\rm Rm}|^{p-1}\phi^{2p-2})^{\frac{p}{p-1}}}{\frac{p}{p-1}}+\frac{(K\rho^{-2}e^{2KT})^p}{p}\right]dV_t\\&\leq A_1+K^pe^{2pKT}p^{-1}\rho^{-2p}{\rm Vol}_{g(t)}\left(B_{g(0)}\left(x_0,\frac{\rho}{\sqrt{K}}\right)\right)\\&\leq U+K^pe^{2pKT}\rho^{-2p}{\rm Vol}_{g(t)}\left(B_{g(0)}\left(x_0,\frac{\rho}{\sqrt{K}}\right)\right)
\end{align*}
Hence
\begin{align*}
    U'\leq \Lambda_1U+\left[CK(K+L^2)K^pe^{2pKT}\rho^{-2p}+CKC^{p}L^2\right]{\rm Vol}_{g(t)}\left(B_{g(0)}\left(x_0,\frac{\rho}{\sqrt{K}}\right)\right),
\end{align*}
where $\Lambda_1 := C(p-1)KL^2+CK(K+L^2)$ is a constant. The Bishop-Gromov volume comparison theorem shows that the inequality
\begin{align*}
    {\rm Vol}_{g(t)}\left(B_{g(0)}\left(x_0,\frac{\rho}{\sqrt{K}}\right)\right)\leq e^{cT}{\rm Vol}_{g(\tau)}\left(B_{g(0)}\left(x_0,\frac{\rho}{\sqrt{K}}\right)\right)
\end{align*}
hols for all $0\leq t\leq \tau\leq T$. consequently, we arrive at
\begin{align*}
    U'\leq \Lambda_1U+\Lambda_2e^{cT}{\rm Vol}_{g(\tau)}\left(B_{g(0)}\left(x_0,\frac{\rho}{\sqrt{K}}\right)\right),
\end{align*}
with $\Lambda_2:=CK(K+L^2)K^pe^{2pKT}\rho^{-2p}+CKC^{p}L^2$. This implies that
\begin{align*}
    \frac{d}{dt}\left(e^{-\Lambda_{1}t}U(t)\right)\leq \Lambda_{2}e^{c(T-t)}{\rm Vol}_{g(\tau)}\left(B_{g(0)}\left(x_0,\frac{\rho}{\sqrt{K}}\right)\right).
\end{align*}
Upon integration over $[0,\tau]$, it yields
\begin{align*}
    U(\tau)\leq e^{\Lambda_1T}\left(U(0)+\Lambda_2{\rm Vol}_{g(\tau)}\left(B_{g(0)}\left(x_0,\frac{\rho}{\sqrt{K}}\right)\right)\right).
\end{align*}
Now we consider \[U(0) = \left(A_1+\frac{CK}{p-1}A_2+KC^{p}\tilde{S}+CKS+C\int_M|{\rm Ric}|^2|{\rm Rm}|^{p-1}\phi^{2p}dV_t\right)_{t=0}\]. We have proved that
\begin{align*}
    A_4\leq A_1+\Lambda_2e^{2pKT}{\rm Vol}_{g(0)}\left(B_{g(0)}\left(x_0,\frac{\rho}{\sqrt{K}}\right)\right).
\end{align*}
According to the definition, it is clear that
\begin{eqnarray*}
    S&=&\int_M|{\rm Rm}|^{p-1}|\nabla u|^2\phi^{2p}dV_t \ \ \leq \ \ L^2A_2,\\
    \widetilde{S} &= &\int_M|\nabla u|^2\phi^{2p}dV_t \ \ \leq \ \ CL^2{\rm Vol}_{g(0)}\left(B_{g(0)}\left(x_0,\frac{\rho}{\sqrt{K}}\right)\right)
\end{eqnarray*}
Applying Young's inequality to $A_2$, we get
\begin{eqnarray*}
    A_2 &=& \int_M|{\rm Rm}|^{p-1}\phi^{2p}dV_t \ \ = \ \ \int_M\left(|{\rm Rm}|^{p-1}\phi^{2p-2}\right)\phi^2dV_t\\
    &\leq&\frac{p-1}{p}\int_M|{\rm Rm}|^p\phi^{2p}dV_t+\frac{1}{p}\int_M\phi^{2p}dV_t\\
    &\leq& A_1+C{\rm Vol}_{g(t)}\left(B_{g(0)}\left(x_0,\frac{\rho}{\sqrt{K}}\right)\right).
\end{eqnarray*}
The obvious estimate
\begin{align*}
    \int_M|{\rm Ric}|^2|{\rm Rm}|^{p-1}\phi^{2p}dV_t \leq K^2A_2
\end{align*}
tells us that
\begin{eqnarray*}
    U(0)&\leq&\left(\frac{CK}{p-1}+CK^2+CKL^2\right)\int_M|{\rm Rm}(g(0))|^{p}\phi^{2p}dV_{g(0)} \\
    &&+ \ \left(\frac{CK}{p-1}+C+CKL^2\right){\rm Vol}_{g(0)}\left(B_{g(0)}\left(x_0,\frac{\rho}{\sqrt{K}}\right)\right)\\
    &=&\Gamma_1e^{-\Lambda_1T} \int_{B_{g(0)}(x_0,\rho/2\sqrt{K})}|{\rm Rm}(g(0))|^p\phi^{2p}dV_{g(0)}\\
    &&+ \ \left(\Gamma_2e^{-\Lambda_1T}-\Lambda_2\right){\rm Vol}_{g(0)}\left(B_{g(0)}\left(x_0,\frac{\rho}{\sqrt{K}}\right)\right),
\end{eqnarray*}
where
$$
\Gamma_1 := e^{\Lambda_1T}\left(\frac{CK}{p-1}+CK^2+CKL^2\right), \ \ \ \Gamma_2:=e^{\Lambda_1T}\left(\frac{CK}{p-1}+C+CKL^2+\Lambda_2\right).
$$
Plug it into the differential inequality and we obtain for $p\geq 2$
\begin{eqnarray*}
&&\int_{B_{g(0)}(x_0,\rho/2\sqrt{K})}|{\rm Rm}({g(t)})|_{g(t)}^{p}dV_{g(t)}\\
&\leq& \Gamma_1\int_{M}|{\rm Rm}({g(0)})|^{p}\phi^{2p}dV_{g(0)}+\Gamma_2{\rm Vol}_{g(0)}\left(B_{g(0)}\left(x_0,\frac{\rho}{\sqrt{K}}\right)\right)\\
&\leq&\Gamma_1\int_{B_{g(0)}(x_0,\rho/2\sqrt{K})}|{\rm Rm}({g(0)})|^{p}dV_{g(0)}+\Gamma_2{\rm Vol}_{g(0)}\left(B_{g(0)}\left(x_0,\frac{\rho}{\sqrt{K}}\right)\right).
\end{eqnarray*}
We finished the proof.
\end{proof}

As it will be needed in the following discussion, We also restate the Theorem \ref{t1.1} to emphasize the power of $p$, which can be easily obtained from $\Gamma_1, \Gamma_2,  \Lambda_1, \Lambda_2$:


\begin{eqnarray*}
    &&-\kern-11pt\int_{B(x_0,\rho/2\sqrt{K})}|{\rm Rm}(g(t))|_{g(t)}^{p}dV_{g(t)}\\
    &:=&\frac{1}{{\rm Vol}\left(B\left(x_0,\frac{\rho}{2\sqrt{K}}\right)\right)}\int_{B(x_0,\rho/2\sqrt{K})}|{\rm Rm}(g(t))|_{g(t)}^{p}dV_{g(t)}\\
    &\leq&\Gamma_1-\kern-11pt\int_{B(x_0,\rho/2\sqrt{K})}|{\rm Rm}(g(0))|^{p}dV_{g(0)}+\Gamma_2\frac{{\rm Vol}_{g(0)}\left(B_{g(0)}\left(x_0,\frac{\rho}{\sqrt{K}}\right)\right)}{{\rm Vol}_{g(0)}\left(B_{g(0)}\left(x_0,\frac{\rho}{2\sqrt{K}}\right)\right)}
    \\
    &\leq& \Gamma_1-\kern-11pt\int_{B(x_0,\rho/2\sqrt{K})}|{\rm Rm}(g(0))|^{p}dV_{g(0)}+C\Gamma_2 e^{C(T+\frac{\rho}{\sqrt{K}})}
    \\
    &\leq& Ce^{C(p-1)}-\kern-11pt\int_{B(x_0,\rho/2\sqrt{K})}|{\rm Rm}(g(0))|^{p}dV_{g(0)}+ Ce^{C(p-1)}K^p\rho^{-2p}
\end{eqnarray*}
where all other constants in it are independent of $p$.

\subsection{Proof of Propositions \ref{p2.1}-\ref{p2.5}}\label{subsection2.2}

In this subsection we give proofs of Proposition \ref{p2.1} -- Proposition \ref{p2.5}.

\begin{proposition}\label{p2.8} We have
\[\frac{d}{dt}A_1\leq B_1+CKB_2+CKA_4+C(K+L^2)A_1+CT_p\]
\end{proposition}

\begin{proof} Compute
\begin{align*}
    \frac{d}{dt}&\left(\int_M|{\rm Rm}|^p\phi^{2p}dV_t\right) = \int_M(\partial_t|{\rm Rm}|^p)\phi^{2p}dV_t +\int_M|{\rm Rm}|^p\phi^{2p}(-{\rm R}+2|\nabla u|^2)dV_t\\&
    =\frac{p}{2}\int_M|{\rm Rm}|^{p-2}[\nabla^2{\rm Ric}\ast {\rm Rm} +{\rm Ric}\ast {\rm Rm}\ast {\rm Rm}\\& \qquad\qquad+{\rm Rm}\ast\nabla^2 u \ast \nabla^2u +{\rm Rm}\ast {\rm Rm}\ast \nabla u\ast \nabla u]\phi^{2p} dV_t \\& \qquad\qquad-\int_M{\rm R}|{\rm Rm}|^p\phi^{2p}dV_t+2\int_M|{\rm Rm}|^p |\nabla u|^2\phi^{2p}dV_t
    \\&\leq C\int_M|{\rm Rm}|^{p-2}(\nabla^2{\rm Ric}\ast {\rm Rm})\phi^{2p}dV_t+CKA_1+CT_p+CL^2A_1
\end{align*}
From (2.5), (2.6) and (2.7) in \cite{KMW2016}, we have:
\begin{align*}
    C\int_M|{\rm Rm}|^{p-2}(\nabla^2 {\rm Ric}\ast {\rm Rm})\phi^{2p}dV_t\leq B_1
    +CKB_2+CKA_4.
\end{align*}
Combine them and we prove the proposition.
\end{proof}

\begin{proposition}\label{p2.9} We have
\begin{align*}B_1\leq CKB_2+C(K+L^2)A_1+CKL^2A_2+CKA_4\\+C_0T_p-\frac{1}{2K}\frac{d}{dt}
\left(\int_M|{\rm Ric}|^2|{\rm Rm}|^{p-1}\phi^{2p}dV_t\right).\end{align*}
\end{proposition}

\begin{proof}From the evolution equation of $|{\rm Ric}|^2$(see \cite{List2005}), we can deduce that
\begin{align*}
    |\nabla {\rm Ric}|^2 =& -\frac{1}{2}\square |{\rm Ric}|^2 +2R_{pijq}R^{pq}R^{ij}-4R_{pijq}R^{ij}\nabla^pu\nabla^qu\\&+4\Delta uR^{ij}\nabla_i\nabla_ju-4R^{ij}\nabla_i\nabla_ku\nabla^k\nabla_ju-4R_{ij}R^j{}_{k}\nabla^iu\nabla u\\&\leq -\frac{1}{2}\square |{\rm Ric}|^2 +CK(L^2+K)|{\rm Rm}|+CK|\nabla^2 u|^2+CK^2L^2,
\end{align*}
in which we used the fact that $|\Delta u|\leq \sqrt{n}|\nabla^2 u|$. Hence we have
\begin{align*}
    B_1&\leq \int_M\Bigg[\frac{1}{2K}(\Delta - \partial_t)|{\rm Ric}|^2+C(L^2+K)|{\rm Rm}| \\& \qquad\qquad+CKL^2+C|\nabla^2 u|^2\Bigg]|{\rm Rm}|^{p-1}\phi^{2p}dV_t\\
    &=\frac{1}{2K}\int_M\left[(\Delta - \partial_t)|{\rm Ric}|^2\right]|{\rm Rm}|^{p-1}\phi^{2p}dV_t\\& \qquad\qquad+C(L^2+K)A_1+CKL^2A_2+CT_p\\
    &=\frac{1}{2K}\int_M(\Delta|{\rm Ric}|^2)|{\rm Rm}|^{p-1}\phi^{2p}dV_t+C(L^2+K)A_1+CKL^2A_2+CT_p\\&
     \qquad\qquad-\frac{1}{2K}\int_M\Bigg[\partial_t(|{\rm Ric}|^2|{\rm Rm}|^{p-1}\phi^{2p}dV_t)\\& \qquad\qquad-|{\rm Ric}|^2(\partial_t|{\rm Rm}|^{p-1})\phi^{2p}dV_t-|{\rm Ric}|^2|{\rm Rm}|^{p-1}\phi^{2p}(-R+2|\nabla u|^2)dV_t\Bigg]\\&
    =-\frac{1}{2K}\left(\int_M\langle \nabla|{\rm Ric}|^2,\nabla|{\rm Rm}|^{p-1}\rangle\phi^{2p}dV_t+\int_M\langle\nabla|{\rm Ric}|^2,\nabla\phi^{2p}\rangle|{\rm Rm}|^{p-1}dV_t\right)\\& \qquad\qquad-\frac{1}{2K}\left(\frac{d}{dt}\int_M|{\rm Ric}|^2|{\rm Rm}|^{p-1}\phi^{2p}dV_t\right)+C(L^2+K)A_1+CKL^2A_2\\& \qquad\qquad+CT_p+\frac{1}{2K}\int_M|{\rm Ric}|^2(\partial_t|{\rm Rm}|^{p-1})\phi^{2p}dV_t
\end{align*}
From the proof of (2.13)-(2.15) in \cite{KMW2016}, we can deduce:
\begin{align*}
    \frac{C}{K}\int_M|{\rm Ric}|^2|{\rm Rm}|^{p-3}\phi^{2p}(\nabla^2{\rm Ric}\ast {\rm Rm})dV_t\leq \frac{1}{5}B_1+CKB_2+CKA_4
\end{align*}
Then we can write:
\begin{align*}
    \frac{1}{2K}\int_M(\partial_t|{\rm Rm}|^{p-1})\phi^{2p}dV_t = \frac{p-1}{4K}\int_M|{\rm Ric}|^2(|{\rm Rm}|^{p-3}\partial_t|{\rm Rm}|^2)\phi^{2p}dV_t\\=\frac{C}{K}\int_M|{\rm Ric}|^2|{\rm Rm}|^{p-3}\phi^{2p}[\nabla^2{\rm Ric}\ast {\rm Rm}+{\rm Ric}\ast {\rm Rm}\ast {\rm Rm}+\\{\rm Rm}\ast \nabla^2u\nabla^2u +{\rm Rm}\ast {\rm Rm}\ast \nabla u\ast\nabla u]dV_t\\\leq \frac{1}{5}B_1+CKB_2+CKA_1+CKL^2A_2+CKA_4+CT_p
\end{align*}
From (2.10) and (2.11) in \cite{KMW2016}, we have:
\begin{align*}
    -\frac{1}{2K}\int_M\langle \nabla|{\rm Ric}|^2,\nabla|{\rm Rm}|^{p-1}\rangle\phi^{2p}dV_t\leq \frac{1}{10}B_1+CKB_2\\
    -\frac{1}{2K}\int_M\langle\nabla|{\rm Ric}|^2,\nabla\phi^{2p}\rangle|{\rm Rm}|^{p-1}dV_t\leq \frac{1}{10}B_1+CKA_4
\end{align*}
Plugging them all together and we arrive at Proposition \ref{p2.2}.
\end{proof}

As already stated in notations that all $C$ are irrelevant constants, while $C_0$ in Proposition \ref{p2.2} is a special constant used latter.

\begin{proposition}\label{p2.10} We have
\[B_2\leq -\frac{1}{p-1}\frac{d}{dt}\left(\int_M|{\rm Rm}|^{p-1}\phi^{2p}dV_t\right)+CA_1+CA_4+CL^2A_2+CT_{p-1}\]
\end{proposition}

\begin{proof} Using the evolution inequality of $|{\rm Rm}|$(see \cite{List2005}), we can obtain:
\begin{align*}
    B_2&\leq \int_M\left[\frac{1}{2}(\Delta-\partial_t)|{\rm Rm}|^2+C|{\rm Rm}|^3+CL^2|{\rm Rm}|^2+C|\nabla^2u|^2|{\rm Rm}|\right]|{\rm Rm}|^{p-3}\phi^{2p}dV_t\\
    &=\frac{1}{2}\int_M(\Delta|{\rm Rm}|^2)|{\rm Rm}|^{p-3}\phi^{2p}dV_t +CA_1+CL^2A_2\\& \qquad\qquad\qquad\qquad\qquad\qquad+T_{p-1}
    -\frac{1}{2}\int_M(\partial_t |{\rm Rm}|^2)|{\rm Rm}|^{p-3}\phi^{2p}dV_t\\&
    \leq C\int_M|\nabla {\rm Rm}||\nabla \phi||{\rm Rm}|^{p-2}\phi^{2p-1}dV_t+CA_1+CL^2A_2\\& \qquad\qquad\qquad\qquad\qquad\qquad+T_{p-1}
    -\frac{1}{2}\int_M(\partial_t |{\rm Rm}|^2)|{\rm Rm}|^{p-3}\phi^{2p}dV_t\\&
    \leq  \frac{1}{2}B_2+CA_4+CA_1+CL^2A_2+T_{p-1} -\frac{1}{2}\int_M(\partial_t |{\rm Rm}|^2)|{\rm Rm}|^{p-3}\phi^{2p}dV_t
\end{align*}
Following the proof of (2.18)-(2.19) in \cite{KMW2016},
\begin{align*}
    -\frac{1}{2}\int_M(\partial_t&|{\rm Rm}|^2)|{\rm Rm}|^{p-3}\phi^{2p}dV_t \\&= -\frac{1}{2}\int_M[\partial_t(|{\rm Rm}|^2|{\rm Rm}|^{p-3}\phi^{2p}dV_t)\\& \qquad\qquad
    -|{\rm Rm}|^2(\partial_t|{\rm Rm}|^{p-3})\phi^{2p}dV_t-|{\rm Rm}|^{p-1}\phi^{2p}\partial_tdV_t]\\&
    =-\frac{1}{2}\partial_tA_2+\frac{p-3}{4}\int_M(\partial_t|{\rm Rm}|^2)|{\rm Rm}|^{p-3}\phi^{2p}dV_t\\& \qquad\qquad-\frac{1}{2}\int_M{\rm R}|{\rm Rm}|^{p-1}\phi^{2p}dV_t+\int_M|{\rm Rm}|^{p-1}|\nabla u|^2\phi^{2p}dV_t.
\end{align*}
Therefore, we can find:
\begin{align*}
    -\frac{1}{2}\int_M(\partial_t|{\rm Rm}|^2)|{\rm Rm}|^{p-3}\phi^{2p}dV_t\leq -\frac{1}{p-1}\partial_t A_2+CA_1+CA_4+CL^2A_2+CT_{p-1}
\end{align*}
In summary we can find
\begin{align*}
    B_2\leq -\frac{1}{p-1}\frac{d}{dt}\left(\int_M|{\rm Rm}|^{p-1}\phi^{2p}dV_t\right)+CA_1+CA_4+CL^2A_2+CT_{p-1}
\end{align*}
and finish the proof.
\end{proof}

\begin{proposition}\label{p2.11} For any $p\geq2$, $T_p$ satisfy the following estimate
\begin{align*}
    T_p\leq-\frac{d}{dt}S-\frac{C}{p-1}&\frac{d}{dt}A_2-\frac{C}{K}\frac{d}{dt}\left(\int_M|{\rm Ric}|^2|{\rm Rm}|^{p-1}\phi^{2p}dV_t\right)\\&+C(K+L^2)A_1+CKL^2A_2+C(K+L^2)A_4+C^{p-1}T_1
\end{align*}
\end{proposition}

\begin{proof} We consider the quantity:
\begin{eqnarray*}
    &&\frac{d}{dt}\left(\int_M|{\rm Rm}|^{p-1}|\nabla u |^2\phi^{2p}dV_t\right)\\
    &=& \int_M(\partial_t|{\rm Rm}|^{p-1})|\nabla u |^2\phi^{2p}dV_t\\
    &&- \ \int_M|{\rm Rm}|^{p-1}|\nabla u|^2({\rm R}-2|\nabla u|^2)\phi^{2p}dV_t\\
    &&+ \ \int_M|{\rm Rm}|^{p-1}(\Delta |\nabla u|^2-2|\nabla^2 u|^2-4|\nabla u|^4)\phi^{2p}dV_t,
\end{eqnarray*}
which infer:
\begin{eqnarray*}
    T_p&=&\int_M|{\rm Rm}|^{p-1}|\nabla^2u|^2\phi^{2p}dV_t\\
    &=&-\frac{1}{2}\frac{d}{dt}\int_M|{\rm Rm}|^{p-1}|\nabla u |^2\phi^{2p}dV_t\\
    &&+ \ \frac{1}{2}\int_M(\partial_t|{\rm Rm}|^{p-1})|\nabla u |^2\phi^{2p}dV_t\\
    &&- \ \frac{1}{2}\int_M|{\rm Rm}|^{p-1}|\nabla u|^2({\rm R}-2|\nabla u|^2)\phi^{2p}dV_t\\
    &&+ \ \int_M|{\rm Rm}|^{p-1}\left(\frac{1}{2}\Delta |\nabla u|^2-2|\nabla u|^4\right)\phi^{2p}dV_t.
\end{eqnarray*}
Using \begin{align}
    \square |\nabla u|^2 = -2|\nabla^2u|^2-4|\nabla u|^4,
\end{align} from \cite{List2005} we yields that ${\rm R}-2|\nabla u|^2 \geq-C$ and then
\begin{align*}
    -\frac{1}{2}\int_M|{\rm Rm}|^{p-1}|\nabla u|^2({\rm R}-2|\nabla u|^2)\phi^{2p}dV_t\leq CS.
\end{align*}
Therefore, we arrive at
\begin{eqnarray*}
    T_p&\leq& -\frac{1}{2}\frac{d}{dt}\int_M|{\rm Rm}|^{p-1}|\nabla u |^2\phi^{2p}dV_t\\
    &&+ \ \frac{1}{2}\int_M|\nabla u|^2(\partial_t|{\rm Rm}|^{p-1})\phi^{2p}dV_t\\
    &&+ \ CS+\frac{1}{2}\int_M|{\rm Rm}|^{p-1}\Delta |\nabla u|^2\phi^{2p}dV_t.
\end{eqnarray*}
Notice that by the evolution equation of $|{\rm Rm}|^2$ (see \cite{List2005})
\begin{align*}
    \frac{1}{2}\int_M|\nabla u|^2&(\partial_t|{\rm Rm}|^{p-1})\phi^{2p}dV_t \\&=\frac{1}{2}\int_M|\nabla u|^2(\nabla^2{\rm Ric}\ast {\rm Rm}+{\rm Ric}\ast {\rm Rm}\ast {\rm Rm}+{\rm Rm}\ast \nabla^2u\ast\nabla^2u \\ &\qquad\qquad+{\rm Rm}\ast {\rm Rm}\ast \nabla u\ast \nabla u)|{\rm Rm}|^{p-3}\phi^{2p}dV_t\\&\leq C\int_M|\nabla u|^2\ast\nabla^2 {\rm Ric}\ast|{\rm Rm}|^{p-2}\phi^{2p}dV_t\\& \qquad\qquad+CL^2A_1+CL^2T_{p-1}+CL^2S\\&= -C\int_M\langle\nabla|\nabla u|^2,\nabla {\rm Ric}\rangle|{\rm Rm}|^{p-2}\phi^{2p}dV_t\\& \qquad\qquad-CL^2\int_M\langle\nabla|{\rm Rm}|^2,\nabla {\rm Ric}\rangle|{\rm Rm}|^{p-4}\phi^{2p}dV_t\\& \qquad\qquad-CL^2\int_M\langle\nabla \phi,\nabla {\rm Ric}\rangle|{\rm Rm}|^{p-2}\phi^{2p-1}dV_t\\& \qquad\qquad+CL^2A_1+CL^2T_{p-1}+CL^2S\\&\leq C\int_M|\nabla^2u||\nabla u||\nabla {\rm Ric}||{\rm Rm}|^{p-2}\phi^{2p}\\& \qquad\qquad+CL^2\int_M|\nabla {\rm Rm}||
    \nabla {\rm Ric}||{\rm Rm}|^{p-3}\phi^{2p}dV_t\\& \qquad\qquad+CL^2\int_M|\nabla \phi||\nabla {\rm Rm}||{\rm Rm}|^{p-2}\phi^{2p-1}dV_t\\& \qquad\qquad+CL^2A_1+CL^2T_{p-1}+CL^2S\\&
   \leq CT_{p-2}+\frac{1}{8C_0}B_1+CL^2B_2+CA_4+CL^2A_1+CL^2T_{p-1}+CL^2S
\end{align*}
Applying integrating by parts, the last term becomes
\begin{align*}
    \int_M|{\rm Rm}|^{p-1}&\Delta |\nabla u|^2\phi^{2p}dV_t \\&= -\int_M\langle \nabla|\nabla u|^2,\nabla|{\rm Rm}|^{p-1}\phi+2p|{\rm Rm}|^{p-1}\nabla\phi\rangle\phi^{2p-1}dV_t\\&\leq 2C\int_M|\nabla^2 u||\nabla u||\nabla {\rm Rm}||{\rm Rm}|^{p-2}\phi^{2p}dV_t\\& \qquad\qquad+2C\int_M|\nabla^2 u||\nabla u||\nabla \phi||{\rm Rm}|^{p-1}\phi^{2p-1}dV_t\\&\leq \frac{1}{8}T_p+8CL^2B_2+\frac{1}{8}T_p+8CL^2A_4
\end{align*}
Plugging them into the inequality of $T_p$, we obtain
\begin{align*}
    T_p\leq-\frac{1}{2}\partial_tS&+ CT_{p-2}+\frac{1}{8C_0}B_1+CL^2A_1\\&+CL^2T_{p-1}+CL^2S+\frac{1}{8}T_p+CL^2B_2+CL^2A_4
\end{align*}
Replacing $B_1$ by using Proposition \ref{p2.9} yields
\begin{align*}
    T_p\leq-\frac{1}{2}\partial_tS&+ CT_{p-2}-\frac{1}{16C_0K}\frac{d}{dt}\left(\int_M|{\rm Ric}|^2|{\rm Rm}|^{p-1}\phi^{2p}dV_t\right)\\&+C(K+L^2)A_1+CKL^2A_2+CL^2T_{p-1}\\&+CL^2S+\frac{1}{4}T_p+CL^2B_2+CL^2A_4
\end{align*}
Using the relationship between $T_k$ (see Lemma \ref{l2.3}), we can write inequalities:
\[CT_{p-2}\leq C\left[\frac{1}{8C}T_p+2(8C)^{\frac{p-3}{2}}T_1\right]\leq\frac{1}{8}T_p+2(8C)^{\frac{p}{2}}T_1\]
to replace $CT_{p-2}$ and we will get:
\begin{align*}
    T_p\leq-\frac{1}{2}\partial_tS&+ \frac{3}{8}T_p+2(8C)^{\frac{p}{2}}T_1+C(K+L^2)A_1\\& +CKL^2A_2-\frac{1}{16C_0K}\frac{d}{dt}\left(\int_M|{\rm Ric}|^2|{\rm Rm}|^{p-1}\phi^{2p}dV_t\right)\\&+CL^2T_{p-1}+CL^2S+CL^2B_2+CL^2A_4
\end{align*}
Replacing $B_2$ by using Proposition \ref{p2.4}, we obtain
\begin{align*}
    T_p\leq-\frac{1}{2}\partial_tS&-\frac{CL^2}{p-1}\frac{d}{dt}A_2-\frac{1}{16C_0K}\frac{d}{dt}\left(\int_M|{\rm Ric}|^2|{\rm Rm}|^{p-1}\phi^{2p}dV_t\right)\\&+ \frac{3}{8}T_p+2(8C)^{\frac{p}{2}}T_1+C(K+L^2)A_1+CKL^2A_2\\&+CL^2T_{p-1}+CL^2S+C(K+L^2)A_4
\end{align*}
Again we can write
\[CL^2T_{p-1}\leq CL^2\left[\frac{1}{8CL^2}T_p+(8CL^2)^{p-2}T_1\right]=\frac{1}{8}T_p+(8CL^2)^{p}T_1\]
Plugging it into the inequality and we finally have
\begin{align*}
    T_p\leq-\frac{1}{2}\partial_tS-\frac{CL^2}{p-1}&\frac{d}{dt}A_2-\frac{1}{16C_0K}\frac{d}{dt}\left(\int_M|{\rm Ric}|^2|{\rm Rm}|^{p-1}\phi^{2p}dV_t\right)\\&+ \frac{1}{2}T_p+C(K+L^2)A_1+CKL^2A_2\\&+C^{p-1}T_1+CL^2S+C(K+L^2)A_4
\end{align*}
which infer:
\begin{align*}
    T_p\leq-&\partial_tS-\frac{CL^2}{p-1}\frac{d}{dt}A_2-\frac{C}{K}\frac{d}{dt}\left(\int_M|{\rm Ric}|^2|{\rm Rm}|^{p-1}\phi^{2p}dV_t\right)\\&+C(K+L^2)A_1+CKL^2A_2+C^{p-1}T_1+CL^2S+C(K+L^2)A_4
\end{align*}
Then we finish the proof.
\end{proof}

\begin{proposition}\label{p2.12} $T_1$ satisfy the following estimate
\begin{align*}
    T_1\leq -\partial_t\widetilde{S}+CL^2{\rm Vol}_{g(t)}(\Omega)
\end{align*}
\end{proposition}

\begin{proof} Consider the quantity: \begin{align*}
    \partial_t\widetilde{S} &= \partial_t\int_M|\nabla u|^2\phi^{2p}dV_t\\&=\int_M(\Delta|\nabla u|^2-2|\nabla^2u|^2-4|\nabla u|^4)\phi^{2p}dV_t+\int_M|\nabla u|^2\phi^{2p}(-{\rm R}+2|\nabla u|^2)dV_t\\&\leq -2T_1+\int_M\Delta|\nabla u|^2\phi^{2p}dV_t+CL^2\int_M\phi^{2p}dV_t\\&\leq -2T_1+2C\int_M|\nabla^2 u||\nabla u||\nabla \phi|\phi^{2p-1}dV_t+CL^2{\rm Vol}_{g(t)}(\Omega)
    \\&\leq -T_1+C\int_M|\nabla u|^2|\nabla \phi|^2\phi^{2p-2}dV_t+CL^2{\rm Vol}_{g(t)}(\Omega)\\&\leq -T_1+CL^2{\rm Vol}_{g(t)}(\Omega)
\end{align*}
\end{proof}

\section{The extension of the Ricci-Harmonic flow}\label{section3}

As \cite{List2008} has proved, the flow can be extended over $T$ if the Riemannian curvature is bounded at each point. First we prove

\begin{lemma}\label{l3.1} There exist constants $C$ such that the following estimate
\begin{align*}
    \square |{\rm Rm}|&\leq C|{\rm Rm}|^2+C|\nabla^2 u|^2+C
\end{align*}
holds.
\end{lemma}

\begin{proof} Using the evolution equation of $|{\rm Rm}|^2$ (see Chapter 2.7 in \cite{List2005}), we obtain:
\begin{align*}
    \square |{\rm Rm}|^2 &= 2|{\rm Rm}|(\partial_t |{\rm Rm}|) - 2|{\rm Rm}|(\Delta |{\rm Rm}|)-2|\nabla|{\rm Rm}||^2 \\
    & = 2|{\rm Rm}|(\square |{\rm Rm}|)-2|\nabla|{\rm Rm}||^2 \\
    &\leq -2|\nabla {\rm Rm}|^2+C|{\rm Rm}|^3+C|{\rm Rm}||\nabla^2u|^2+C|\nabla u|^2|{\rm Rm}|^2
\end{align*}
From $|\nabla {\rm Rm}|\geq |\nabla |{\rm Rm}||$ and assumption (2), we can get
\begin{align*}
    \square |{\rm Rm}| &\leq C|{\rm Rm}|^2 + C|\nabla^2 u|^2 + CL^2|{\rm Rm}|\\
    &\leq C|{\rm Rm}|^2  + C|\nabla^2 u|^2 + CL^2(|{\rm Rm}|^2+1)\\
    &= C|{\rm Rm}|^2  + C|\nabla^2 u|^2 + C
\end{align*}
which gives the desired estimate.
\end{proof}

\par Now we prove Theorem \ref{t1.4}.

\begin{theorem}\label{t3.2} Let $(g(t),u(t))$ be a smooth solution to the Ricci-harmonic flow on $M\times[0,T)$ with $T<\infty$, where $M$ is a complete $n$-dimensional manifold. If $(M,g(0))$ is complete and:
    \begin{align*}
        \underset{\quad M\quad }{\sup}|{\rm Rm}(g(0))|_{g(0)}<\infty, \ \ \  \sup_{M\times[0,T)}|{\rm Ric}(g(t))|_{g(t)}<\infty,
    \end{align*}
    then $|{\rm Rm}|$ is locally bounded and $g(t)$ extends smoothly to a complete solution on $[0,T+\epsilon)$ for some constants $\epsilon>0$.
\end{theorem}

\begin{proof} According to Remark \ref{r1.2}, we can denote
\begin{align*}
        K:= \sup_{M\times[0,T)}|{\rm Ric}|(x,t)<\infty, \ \ \ L:= \sup_{M\times[0,T)}|\nabla u|(x,t)<\infty.
\end{align*}
According to Lemma \ref{l3.1}, we can pick a constant $C_m\geq 2$ that is sufficiently
large so that
\begin{align*}
    \square |{\rm Rm}|\leq C_m(|{\rm Rm}|^2+2|\nabla^2 u|^2+1)
\end{align*}
Plugging it with evolution equation (2.1) we can find
\begin{align*}
    (\partial_t-\Delta )(|{\rm Rm}|+C_m|\nabla u|^2+1)&=(\partial_t-\Delta )(|{\rm Rm}|+C_m|\nabla^2 u|^2)\\
    &=C_m(|{\rm Rm}|^2-4|\nabla u|^4+1)\\
    &\leq C_m(|{\rm Rm}|^2+C_m^2|\nabla u|^4+1)\\
    &\leq C_m(|{\rm Rm}|+C_m|\nabla u|^2+1)^2
\end{align*}
On the other hand,
\begin{align*}
    \int_\Omega(|{\rm Rm}|+C_m|\nabla u|^2+1)^pdV_{g(t)}\leq 3^{p-1}\int_\Omega(|{\rm Rm}|^p+C_m^p|\nabla u|^{2p}+1)dV_{g(t)}\\
    \leq 3^{p-1}\int_\Omega|{\rm Rm}|^pdV_{g(t)} +3^{p-1}(C_m^pL^{2p}+1){\rm Vol}_{g(t)}(\Omega)
\end{align*}
Define
\begin{align*}
    {\Phi} := |{\rm Rm}|+C_m|\nabla u|^2 +1
\end{align*}
and then the above propositions gives
\begin{align*}
    \left(-\kern-11pt\int_\Omega {\Phi}^p dV_{g(t)}\right)^{\frac{1}{p}}&\leq \left(3^{p-1}-\kern-11pt\int_\Omega|{\rm Rm}|^pdV_{g(t)}+3^{p-1}(C_m^pL^{2p}+1)\right)^{\frac{1}{p}}\\
    &\leq 3\left(-\kern-11pt\int_\Omega|{\rm Rm}|^pdV_{g(t)}\right)^{\frac{1}{p}}+3C_mL^2+3\\
    &\leq 3\left[Ce^{C(p-1)}(\Lambda+K^p\rho^{-2p})\right]^{\frac{1}{p}}+3C_mL^2+3\\
    &\leq C(1+\Lambda)+3K\rho^{-2}+ 3C_mL^2+3
    \\&:=C_n,
\end{align*}
which is a constant independent of $p$. We also have
\begin{align*}
    (\partial_t - \Delta){\Phi}\leq C_m {\Phi}^2.
\end{align*}
The progress to give uniform bound from $L^p$ estimate is an essentially routine applying De Giorgi-Nash-Moser iteration presented in Lemma 19.1 of \cite{LP2012}. We write $f = u = {\Phi}$ and the above inequality shows that
\begin{align*}
    \partial_t u \leq \Delta u +Cfu
\end{align*}
weakly on $M\times [0,T]$. It is equivalent to say that for fixed $a\geq 1$
\begin{align}
    -\int_M \varphi^2u^{2a-1}\Delta u\!\ dV_{g(t)}+\frac{1}{2a}\int_M\varphi^2\partial_t(u^{2a})dV_{g(t)}\leq C\int_M \varphi^2u^{2a}f\!\ dV_{g(t)}\label{3.1}
\end{align}
for any $t\in [0,T]$ and non-negative Lipschitz function $\varphi$ whose support is compactly contained in $B_{g(0)}(x_0,\rho/2\sqrt{K})$. Integrate by part and notice that $a\geq 1$, we obtain
\begin{align*}
     -\int_M \varphi^2u^{2a-1}&\Delta u\!\ dV_{g(t)} \\
     &= 2\int_M\varphi u^{2a-1}\langle\nabla u,\nabla \varphi\rangle dV_{g(t)}+(2a-1)\int_M\varphi^2u^{2a-2}|\nabla u|^2dV_{g(t)}\\&
     \geq \frac{1}{a}\int_M2a\varphi u^{2a-1}\langle\nabla u,\nabla \varphi\rangle dV_{g(t)}+\frac{1}{a}\int_Ma^2\varphi^2u^{2a-2}|\nabla u|^2dV_{g(t)}\\&=\frac{1}{a}\int_M |\nabla(\varphi u^a)|^2dV_{g(t)} -\frac{1}{a}\int_M |\nabla \varphi|^2u^{2a}dV_{g(t)}
\end{align*}
For Ricci-Harmonic flow, we have $\partial_t dV_{g(t)} = (-R+2|\nabla u|^2)dV_{g(t)}$, and furthermore
\begin{align*}
    \left|R-2|\nabla u|^2\right|\leq |R|+2|\nabla u|^2
    \leq C\left(|{\rm Rm}|+C_m|\nabla u|^2+1\right)=C{\Phi}=Cf,
\end{align*}
we then arrive at
\begin{align*}
    \int_M\varphi^2\partial_t(u^{2a})dV_{g(t)} &= \frac{d}{dt}
    \left(\int_M\varphi^2u^{2a}dV_{g(t)}\right)- \int_M\varphi^2u^{2a}(R-2|\nabla  u|^2)dV_{g(t)}\\&\geq  \frac{d}{dt}\left(\int_M\varphi^2u^{2a}dV_{g(t)}\right)- C\int_M\varphi^2u^{2a}f\!\ dV_{g(t)}.
\end{align*}
Plugging the above two inequalities into (\ref{3.1}) implies
\begin{align*}
    \int_M |\nabla(\varphi u^a)|^2dV_{g(t)}+& \frac{1}{2}\frac{d}{dt}\left(\int_M\varphi^2u^{2a}dV_{g(t)}\right)
    \\&\leq Ca\int_M\varphi^2u^{2a}fdV_{g(t)}+\int_M |\nabla \varphi|^2u^{2a}dV_{g(t)}.
\end{align*}
Following (3.6)-(3.11) of \cite{KMW2016} for the rest of the steps with $B = B_{g(0)}(x_0,\rho/2\sqrt{K})$, we obtain the following inequality
\begin{align*}
    \sup_{B_{g(0)}(x_0,\frac{\rho}{4\sqrt{K}})\times [\frac{T}{2},T]} u\leq Ce^{C(T+\frac{\rho}{\sqrt K})}\left(A^\alpha+\left(\left(\frac{\rho}{\sqrt K}\right)^{-2}+T^{-1}\right)\right)^{\frac{2\mu-1}{p(\mu-1)}}A,
\end{align*}
where $\alpha = \frac{p(\mu-1)}{\mu(p-1)-p}$ and $\mu = \mu(n)\leq \frac{n}{n-2}$ is given by the Sobolev inequality (see \cite{KMW2016}). $A$ is the average $L^p$ estimate of $f$, i.e.
\begin{align*}
    A:= \sup_{t\in [0,T]}\left(-\kern-11pt\int_Bf^p(t)dV_0\right)^{\frac{1}{p}}
\end{align*}
Apply the following result back to ${\Phi}$ and we get the local uniform bound for ${\Phi}$ near $T$:
\begin{align*}
    \sup_{B_{g(0)}(x_0,\frac{\rho}{4\sqrt{K}})\times [\frac{T}{2},T]}{\Phi} \leq Ce^{C(T+\frac{\rho}{\sqrt K})}\left(1+C_n^{\alpha'}+\left(\frac{K}{\rho^2}+T^{-1}\right)^{\beta'}\right),
\end{align*}
where constants $\alpha',\beta'$ only depend on $n$ and other constants may depend on $n,K$, $L$, $\rho,\Lambda,C_m$ but not $p$. Finally, since:
\begin{align*}
    \underset{t\to T}{\rm lim} |{\rm Rm}| \leq \underset{t\to T}{\rm lim} {\Phi} < \infty
\end{align*}
satisfied and by the Theorem 6.22 of \cite{List2005}, we immediately yield that the the Ricci-Harmonic flow can be smoothly extended past $T$.
\end{proof}



\end{document}